\documentclass[final,3p,times]{elsarticle}
\usepackage{amsfonts}

 \usepackage{graphics}
 \usepackage{graphicx}
 \usepackage{epsfig}

\usepackage{amssymb}
\usepackage{amsmath,latexsym}
\usepackage{mathrsfs}
\usepackage{bm}
\usepackage{amssymb}
\usepackage{hyperref}
\newtheorem{theorem}{Theorem}[section]

\newtheorem{lem}{Lemma}[section]

\newtheorem{assu}{Assumption}[section]

\numberwithin{equation}{section}

\numberwithin{table}{section}

\numberwithin{figure}{section}

\journal{DCDS}

\begin{document}

\begin{frontmatter}



\title{A General Stochastic Maximum Principle For Optimal Control Of Stochastic Systems Driven By Multidimensional
Teugel's Martingales\tnoteref{t1}} \tnotetext[t1]{This
research was supported by the National Natural Science Foundation of China (Program No.11171215).}

\author[sjtu]{Jianzhong Lin\corref{cor1}}

\ead{jzlin@sjtu.edu.cn} \cortext[cor1]{Corresponding author.}

\address[sjtu]{Department of Mathematics, Shanghai Jiaotong University, Shanghai 200240, China}

\begin{abstract}
A necessary maximum principle is proved for optimal controls of
stochastic systems driven by multidimensional Teugel's martingales. The multidimensional Teugel's martingales are
constructed by orthogonalizing the multidimensional L\'{e}vy processes. The control domain need not be convex,
and the control is allowed to enter into the terms of Teugel's martingales.
\end{abstract}

\begin{keyword}
Stochastic optimal control\sep Maximum principle\sep
Backward stochastic differential equation\sep L\'{e}vy processes\sep Teugel's martingals.

{\bf MSC Subject Classification 2000}:  93E20\sep 60H10\sep 60J75\sep 60G44\sep 49K45
\end{keyword}

\end{frontmatter}


\section{Introduction}\label{sec1}

The stochastic maximum principle is one of the central topics in the stochastic optimal control theory. In the past four decades, a variety of results have been obtained on optimal stochastic control problems.(cf. for example, [1], [3], [5], [12], [14]-[17], [26], [31]). Two major advances in these works are worth mentioning. One is the definition of the adjoint processes and its characterization by It\^{o}-type equations. This was contributed by Kushner [17] and Bismut [5], and summarized by Bensoussan [3] via functional analysis methods. Another advance is the idea of second-order variation in calculating the variation of the cost functional
caused by the spike variation of the given optimal control. This was motivated by the study of the nonconvex optimal stochastic control of diffusion processes with the control entering into the diffusion term, and was developed by Peng [26]. On nonconvex controls of diffusion processes, we refer the reader to Kushner [17], Haussmann [14], Bensoussan [3], Hu [15], Hu and Peng [16], Peng [26] and Yong and Zhou [35].

It is well known that jump-diffusion process is an important class of processes for describing financial data. The stochastic maximum principle of jump-diffusion processes, where the control is unallowed into the jump terms, was considered by Boel [6], Boel and Varaiya [7], Rishel [28], Davis and Elliott [9] and Situ [31]. The further profound problem, where the control enters
into the diffusion and jump terms and also some state constrains are imposed, was completely solved by Tang and Li [34] by applying the idea of second-order variation. On the convex controls of jump-diffusion, we refer the reader to Cadenillas[8], Framstad, Okesendal and Sulem [13], Shi and Wu [30].

The L\'{e}vy process (refers to Bertoin [4], Sato [29]) is a more general class of discontinuous processes than jump-diffusion processes. Nualart and Schoutens [22] obtained some interesting results. They introduce the power jump processes and the related
Teugel's martingales. Furthermore, they give a chaotic and predictable representation for a one-dimensional L\'{e}vy process, in terms of these orthogonalized Teugels martingales. Thus the martingale representation theorem for L\'{e}vy process satisfying some exponential moment condition was a consequence of the chaotic representation. Nualart and Schoutens [23] established the existence and uniqueness of solutions for BSDE driven by a one-dimensional L\'{e}vy process of the kind considered in Nualart and Schoutens [22]. Further progresses on the subject were subsequently given by Bahlali, Eddahbi and Essaky [2], Ren[27], Lin[20]. Based on these Results, a stochastic linear-quadratic problem with L\'{e}vy processes was considered by Mitsui and Tabata [24],Tang and Wu [32]. The stochastic maximum principle, where the control enters into the diffusion and jump terms and also control domain is convex, was given by Meng and Tang [21], Tang and Zhang [33].

Recently, A chaotic and predictable representation theorem associated with multidimensional L\'{e}vy processes was obtained by Lin [19]. This extends the setting in Nualart and Schoutens [22] into the multidimensional L\'{e}vy processes. Furthermore, The existence and uniqueness of solutions for BSDEs driven by multidimensional Teugel's martingales, which are constructed by orthogonalizing the multidimensional L\'{e}vy processes, was proved by Lin [20]. According to these results and following the research line of the paper in Peng [26] and Tang and Li [34], this paper discusses the general stochastic maximum principle where the control systems are driven by the multidimensional Teugel's martingales. It is worth emphasizing that there are three main differences in our setting compared with Mitsui and Tabata [24],Tang and Wu [32],Meng and Tang [21] and Tang and Zhang [33]. First, in our paper, the each component in stochastic system is driven by a Teugel's martingale which is generated by the multidimensional L\'{e}vy processes, while the each component in stochastic system in [21], [24], [32] and [33] is driven by a Teugel's martingale which is generated by one component of multidimensional L\'{e}vy processes. Secondly, in our paper, the control domain need not be convex, while that in Meng and Tang [21], Tang and Zhang [33] is convex and therefore the second-order variation technique is unnecessary. Finally, the terminal state in our case is constrained while is not in Meng and Tang [21], Tang and Zhang [33].

The paper is organized as follows. Section 2 contains an introduction on chaotic and predictable representation theorem associated with multidimensional L\'{e}vy processes and BSDEs driven by multidimensional Teugel's martingales. In Section 3, we give the
statement of the problem, our main assumptions and some preliminary lemmas about the first- and second-order variational equation and variational inequality which will be used in the sequel. In Section 4, we derive the first- and second-order adjoint equations, and finally prove the necessary maximum principle. The conclusions are drawn in Section 5.

\section{BSDE driven by multidimensional Teugel's martingales}\label{sec2}

A $\mathbb{R}^n$-valued stochastic process
$X=\{X(t)=(X_1(t),X_2(t),\cdots,X_n(t))',t\geq 0\}$ defined in complete probability space $(\Omega,\mathscr{F},\mathbb{P})$ is
called \textsl{L\'{e}vy process} if $X$ has stationary and independent increments and $X(0)=\bm{0}$. A L\'{e}vy process
possesses a c\`{a}dl\`{a}g modification and we will always assume that we are using this c\`{a}dl\`{a}g version. If we let
$\mathscr{F}_t=\mathscr{G}_t\vee\mathscr{N}$, where $\mathscr{G}_t=\sigma\{X(s),0\leq s\leq t\}$ is the natural filtration of $X$, and $\mathscr{N}$ are the $\mathbb{P}-$null sets of $\mathscr{F}$, then $\{\mathscr{F}_t,t\geq 0\}$ is a right continuous family of $\sigma-$fields. We assume that $\mathscr{F}$ is generated by $X$. For an up-to-date and comprehensive account of L\'{e}vy processes we refer the reader to Bertoin [4] and Sato [29].

Let $X$ be a L\'{e}vy process and denote by
\begin{eqnarray}
X(t-)=\lim\limits_{s\rightarrow t,s<t}X(s),\quad t>0 ,\nonumber
\end{eqnarray}
the left limit process and by $\triangle X(t)=X(t)-X(t-)$ the jump
size at time $t$. It is known that the law of $X(t)$ is
\textsl{infinitely divisible} with characteristic function of the
form
\begin{eqnarray}
E\left[exp(i\bm{\theta}\cdot
X(t))\right]=\left(\phi(\bm{\theta})\right)^t,\quad
\bm{\theta}=(\theta_1,\theta_2,\cdots,\theta_n)\in
\mathbb{R}^n\nonumber
\end{eqnarray}
where $\phi(\bm{\theta})$ is the characteristic function of
$\bm{X}(1)$. The function $\psi(\bm{\theta})=log\phi(\bm{\theta})$
is called the \textsl{characteristic exponent} and it satisfies the
following famous L\'{e}vy-Khintchine formula (Bertoin, [4]):
\begin{eqnarray}
\psi(\bm{\theta})=-\frac{1}{2}\bm{\theta}\cdot
\Sigma\bm{\theta}+\textrm{i}\bm{a}\cdot
\bm{\theta}+\int_{\mathbb{R}^n}\left(exp(i\bm{\theta}\cdot
\bm{x})-1-\textrm{i}\bm{\theta}\cdot \bm{x} 1_{|\bm{x}|\leq
1}\right)\nu(d\bm{x}).\nonumber
\end{eqnarray}
where $\bm{a},\bm{x}\in \mathbb{R}^n$, $\Sigma$ is a symmetric
nonnegative-definite $n\times n$ matrix, and $\nu$ is a measure on
$\mathbb{R}^n\backslash\{o\}$ with $\int(\|\bm{x}\|^2\wedge
1)\nu(d\bm{x})<\infty$. The measure $\nu$ is called the
\textsl{L\'{e}vy measure} of $X$.

Throughout this paper, we will use the standard multi-index
notation. We denote by $\mathbb{N}_0$ the set of nonnegative
integers. A multi-index is usually denoted by $\bm{p}$,
$\bm{p}=(p_1,p_2,\cdots,p_n)\in\mathbb{N}_0^n$. Whenever $\bm{p}$
appears with subscript or superscript, it means a multi-index. In
this spirit, for example, for $\bm{x}=(x_1,\cdots,x_n)$, a monomial
in variables $x_1,\cdots,x_n$ is denoted by
$\bm{x}^{\bm{p}}=x_1^{p_1}\cdots x_n^{p_n}$.  In addition, we also
define $\bm{p}!=p_1!\cdots p_n!$ and $|\bm{p}|=p_1+\cdots+p_n$; and
if $\bm{p}$, $\bm{q}\in\mathbb{N}_0^n$, then we define
$\delta_{\bm{p},\bm{q}}=\delta_{\bm{p}_1,\bm{q}_1}\cdots\delta_{\bm{p}_n,\bm{q}_n}$.

In the remaining of the paper, we will suppose that
\begin{assu}
the L\'{e}vy measure satisfies for some
$\varepsilon>0$, and $\lambda>0$,
\begin{eqnarray}
\int_{|\bm{x}|\geq \epsilon}exp(\lambda
\|\bm{x}\|)\nu(d\bm{x})<\infty. \nonumber
\end{eqnarray}
\end{assu}

This implies that
\begin{eqnarray}
\int \bm{x}^{\bm{p}}\nu(d\bm{x})<\infty. \quad |\bm{p}|\geq
2\nonumber
\end{eqnarray}
and that the characteristic function $E\left[exp(i\bm{\theta}\cdot
X(t))\right]$ is analytic in a neighborhood of origin $\bm{o}$. As a
consequence, $X(t)$ has moments of all orders and the polynomials
are dense in $L^2(\mathbb{R}^n,\mathbb{P}\circ X(t)^{-1})$ for all
$t>0$.

Fix a time interval $[0,T]$ and set
$L_T^2=L^2(\Omega,\mathscr{F}_T,\mathbb{P})$. We will denote by
$\mathscr{P}$ the predictable sub-$\sigma$-field of
$\mathscr{F}_T\otimes\mathscr{B}_{[0,T]}$. First we introduce some
notation:
\begin{itemize}
\item [$\bullet$]: Let $H_T^2$ denote the space of square integrable
and $\mathscr{F}_t-$progressively one-dimensional measurable
processes $\phi=\{\phi(t),t\in[0,T]\}$ such that
\begin{eqnarray}
\|\phi\|^2=\mathbb{E}\left[\int_0^T\|\phi(t)\|^2dt\right]<\infty.\nonumber
\end{eqnarray}
\item [$\bullet$]: $M_T^2$ will denote the subspace of $H_T^2$
formed by predictable processes.
\item [$\bullet$]: $(H_T^2(l^2))^m$ and $(M_T^2(l^2))^m$ are the
corresponding spaces of $m-$dimensional $l^2-$valued processes
equipped with the norm
\begin{eqnarray}
\|\bm{\phi}_k\|^2_{l^2}&=&\mathbb{E}\left[\int_0^T\sum\limits_{d=1}^{\infty}\sum\limits_{\bm{p}\in\mathbb{N}_d^n}
|\phi_k^{\bm{p}}|^2\right]\qquad k=1,2,\cdots,m,\nonumber\\
\|\bm{\phi}\|_{(l^2)^m}^2&=&\sum\limits_{k=1}^m\|\bm{\phi}_k(t)\|^2_{l^2},\nonumber
\end{eqnarray}
where $\bm{\phi}=(\bm{\phi}_1,\bm{\phi}_2,\cdots,\bm{\phi}_m)'$,
$\bm{\phi}_k=\{\phi_k^{\bm{p}}:\bm{p}\in\mathbb{N}_0^n\}$,
$k=1,2,\cdots,m$ and
$\mathbb{N}_d^n\stackrel{\rm{def}}{=}\{\bm{p}\in\mathbb{N}_0^n:|\bm{p}|=d\}$.
\item [$\bullet$]: Set $\mathcal{H}_T^2=H_T^2\times(M_T^2(l^2))^m$.
\end{itemize}

Following Lin [19] we introduce power jump monomial processes of
the form
\begin{eqnarray}
X(t)^{(p_1,\cdots,p_n)}\stackrel{\rm{def}}{=}\sum\limits_{0<s\leq
t}(\triangle X_1(s))^{p_1}\cdots(\triangle X_n(s))^{p_n},\nonumber
\end{eqnarray}
The number $|\bm{p}|$ is called the total degree of $X(t)^{\bm{p}}$.
Furthermore define
\begin{eqnarray}
Y(t)^{(p_1,\cdots,p_n)}\stackrel{\rm{def}}{=}X(t)^{(p_1,\cdots,p_n)}-\mathbb{E}[X(t)^{(p_1,\cdots,p_n)}]
=X(t)^{(p_1,\cdots,p_n)}-m_{\bm{p}}t ,\nonumber
\end{eqnarray}
the compensated power jump process of multi-index
$\bm{p}=(p_1,p_2,\cdots,p_n)$ where
$m_{\bm{p}}=\int\prod\limits_{i=1}^n x_i^{p_i}\nu(d\bm{x})$. Under
hypothesis 1, $Y(t)^{(p_1,\cdots,p_n)}$ is a normal martingale,
since for an integrable L\'{e}vy process $Z$, the process
$\{Z_t-E[Z_t], t\geq 0\}$ is a martingale. We call
$Y(t)^{(p_1,\cdots,p_n)}$ the \textsl{Teugels martingale monomial}
of multi-index $(p_1,\cdots,p_n)$.

We can apply the standard Gram-Schmidt process with the graded
lexicographical order to generate a biorthogonal basis
$\{H^{\bm{p}},\bm{p}\in\mathbb{N}^n\}$, such that each
$H^{\bm{p}}(|\bm{p}|=d)$ is a linear combination of the
$Y^{\bm{q}}$, with $|\bm{q}|\leq |\bm{p}|$ and the leading
coefficient equal to $1$. We set
\begin{eqnarray}
H^{\bm{p}}&=&Y^{\bm{p}}+\sum\limits_{\bm{q}\prec\bm{p},|\bm{q}|=|\bm{p}|}c_{\bm{q}}Y^{\bm{q}}
+\sum\limits_{k=1}^{|\bm{p}|-1}\sum\limits_{|\bm{q}|=k}c_{\bm{q}}Y^{\bm{q}},\nonumber
\end{eqnarray}
where $\bm{p}=\{p_1,\cdots,p_n\}$, $\bm{q}=\{q_1,\cdots,q_n\}$ and
$\prec$ represent the relation of graded lexicographical order
between two multi-indexes. Some details about the technique and
theory of orthogonal polynomials of several variables refer to Dunkl
and Xu [11].

Set
\begin{eqnarray}
\textsl{p}(\bm{x})^{\bm{p}}&=&\bm{x}^{\bm{p}}+\sum\limits_{\bm{q}\prec\bm{p},|\bm{q}|=|\bm{p}|}c_{\bm{q}}\bm{x}^{\bm{q}}
+\sum\limits_{k=1}^{|\bm{p}|-1}\sum\limits_{|\bm{q}|=k}c_{\bm{q}}\bm{x}^{\bm{q}},\nonumber\\
\tilde{\textsl{p}}(\bm{x})^{\bm{p}}&=&\bm{x}^{\bm{p}}+\sum\limits_{\bm{q}\prec\bm{p},|\bm{q}|=|\bm{p}|}c_{\bm{q}}\bm{x}^{\bm{q}}
+\sum\limits_{k=2}^{|\bm{p}|-1}\sum\limits_{|\bm{q}|=k}c_{\bm{q}}\bm{x}^{\bm{q}},\nonumber
\end{eqnarray}

Set
\begin{eqnarray}
H^{\bm{p}}(t)&=&\sum\limits_{0<s\leq t}\left((\triangle
X_1)^{p_1}\cdots(\triangle
X_n)^{p_n}+\sum\limits_{\bm{q}\prec\bm{p},|\bm{q}|=|\bm{p}|}c_{\bm{q}}(\triangle
X_1)^{q_1}\cdots(\triangle X_n)^{q_n}\right.\nonumber\\
&&\left.+\sum\limits_{k=1}^{|\bm{p}|-1}\sum\limits_{|\bm{q}|=k}c_{\bm{q}}(\triangle
X_1)^{q_1}\cdots(\triangle
X_n)^{q_n}\right),\nonumber\\
&&-t\mathbb{E}\left[X^{\bm{p}}(1)+\sum\limits_{\bm{q}\prec\bm{p},|\bm{q}|=|\bm{p}|}c_{\bm{q}}X^{\bm{q}}(1)
+\sum\limits_{k=1}^{|\bm{p}|-1}\sum\limits_{|\bm{q}|=k}c_{\bm{q}}X^{\bm{q}}(1)\right]\nonumber\\
&=&\left(c_{\bm{e}_1}X_1(1)+\cdots+c_{\bm{e}_n}X_n(1)\right)+\sum\limits_{0<s\leq
t}\tilde{\textsl{p}}(\triangle
X(s))\nonumber\\
&&-t\mathbb{E}\left[\sum\limits_{0<s\leq
t}\tilde{\textsl{p}}(\triangle
X(s))\right]-t\mathbb{E}\left[c_{\bm{e}_1}X_1(1)+\cdots+c_{\bm{e}_n}X_n(1)\right].\nonumber
\end{eqnarray}

Specially we have
\begin{eqnarray}
H^{\bm{e}_1}(t)&=&c_{\bm{e}_1}(1)(X_1(t)-t\mathbb{E}(X_1(1))),\nonumber\\
H^{\bm{e}_2}(t)&=&c_{\bm{e}_2}(2)(X_2(t)-t\mathbb{E}(X_2(1)))+c_{\bm{e}_1}(2)(X_1(t)-t\mathbb{E}(X_1(1))),\nonumber\\
&\vdots&\\
H^{\bm{e}_n}(t)&=&c_{\bm{e}_n}(n)(X_n(t)-t\mathbb{E}(X_n(1)))+c_{\bm{e}_{n-1}}(n)(X_{n-1}(t)-t\mathbb{E}(X_{n-1}(1)))\nonumber\\
&&+\cdots+c_{\bm{e}_1}(n)(X_1(t)-t\mathbb{E}(X_1(1))).\nonumber
\end{eqnarray}

The main tool in the theory of BSDEs is the martingale
representation theorem (cf. Pardoux and Peng [25]). Nualart and Schoutens [22] had proved the representation theorem associated with one-dimensional L\'{e}vy process, furthermore Nualart and Schoutens [23] had established the existence and uniqueness of solutions for BSDE driven by a one-dimensional Teugel's martingale generated by the L\'{e}vy process. The main results in Lin [19] is the Predictable Representation Property (PRP) associated multidimensional L\'{e}vy processes:
\begin{lem}
 Every random variable $F$ in
$L^2(\Omega,\mathscr{F})$ has a representation of the form
\begin{eqnarray}
\begin{array}{rl}
F=&\mathbb{E}(F)+\sum\limits_{d=1}^{\infty}\sum\limits_{\bm{p}\in\mathbb{N}_d^n}
\int_0^T \Phi^{\bm{p}}(s)dH^{\bm{p}}(s)
\end{array}\nonumber
\end{eqnarray}
where $\Phi^{\bm{p}}(s)$ is predictable. This result is an extended
version for the corresponding Theorem in Nualart and Schouten
[22].
\end{lem}

Taking into account the results and notation presented in the
previous section, it seems natural to consider the BSDEs with the following form
\begin{eqnarray}
-d\bm{Y}(t)&=&\bm{f}(t,\bm{Y}(t-),\bm{Z}(t))dt-\sum\limits_{d=1}^{\infty}\sum\limits_{\bm{p}\in\mathbb{N}_d^n}
\bm{z}^{\bm{p}}(s)dH^{\bm{p}}(s), \quad \bm{Y}(T)=\bm{\xi},
\end{eqnarray}
where
\begin{itemize}
\item [$\bullet$]:$\bm{Y}(t)=(Y_1(t),Y_2(t),\cdots,Y_m(t))'$.
\item [$\bullet$]: $\bm{Z}(t)=\{\bm{z}^{\bm{p}}(t)\}_{\bm{p}\in\mathbb{N}_0^n}$, each
component $\bm{z}^{\bm{p}}(t)=(z_1^{\bm{p}},\cdots,z_m^{\bm{p}})'$
is a $m-$variables $\mathscr{F}_t$ predictable function;
\item [$\bullet$] $\bm{f}=(f_1,f_2,\cdots,f_m)':\Omega\times [0,T]\times\mathbb{R}^m\times
  \left(M_T^2(l^2)\right)^m\rightarrow\mathbb{R}^{m}$ is a measurable $m-$dimensional vector function such
  that $\bm{f}(\cdot,\bm{0},\bm{0})\in (H_T^2)^m$.
  \item [$\bullet$] $\bm{f}$ is uniformly Lipschitz in the first two
  components, i.e., there exists $C_k>0$, $k=1,2,\cdots,m$, such that $dt\otimes d\mathbb{P}$
  a.s., for all $(\bm{y}_1,\bm{z}_1)$ and $(\bm{y}_2,\bm{z}_2)$ in $\mathbb{R}^m\times(\bm{l}^2)^m$
  \begin{eqnarray}
   \left|f_k(t,\bm{y}_1,\bm{z}_1)-f_k(t,\bm{y}_2,\bm{z}_2)\right|\leq
   C_k\left(\|\bm{y}_1-\bm{y}_2\|_2+\|\bm{z}_1-\bm{z}_2\|_{(l^2)^m}\right),\qquad k=1,2,\cdots,m.\nonumber
  \end{eqnarray}
  \item [$\bullet$] $\bm{\xi}\in L_T^2(\Omega,\mathbb{P})$.
\end{itemize}

If $(\bm{f},\bm{\xi})$ satisfies the above assumptions, the pair
$(\bm{f},\bm{\xi})$ is said to be \textbf{standard data} for BSDE. A
solution of the BSDE is a pair of processes,
$\{(\bm{Y}(t),\bm{Z}(t)),0\leq t\leq T\}\in H_T^2\times
\left(M_T^2(l^2)\right)^m$ such that the following relation holds
for all $t\in[0,T]$:
\begin{eqnarray}
\bm{Y}(t)=\bm{\xi}+\int_t^T\bm{f}(s,\bm{Y}(s-),\bm{Z}(s))ds-\sum\limits_{d=1}^{\infty}\sum\limits_{\bm{p}\in\mathbb{N}_d^n}
\int_t^T\bm{z}^{\bm{p}}(s)dH^{\bm{p}}(s).
\end{eqnarray}

A key-result concerning the existence uniqueness of solution of BSDEs (2.2) is given by Lin [20]:
\begin{lem}
Given standard data $(\bm{f},\bm{\xi})$, there exists a unique
solution $(\bm{Y},\bm{Z})$ which solves the BSDE (2.3)
\end{lem}

\section{Notations and preliminary lemmas}
\label{sec3}
Consider the following stochastic control system:
\begin{eqnarray}
dx(t)&=&g(x(t-),v(t))dt+\sum\limits_{d=1}^{\infty}\sum\limits_{\bm{p}\in\mathbb{N}_d^n}
\gamma^{\bm{p}}(x(t-),v(t))dH^{\bm{p}}(t),\nonumber\\
x(0)&=&x_0.
\end{eqnarray}
Here and hereafter
\begin{eqnarray}
g(x,v)&:&\mathbb{R}^m\times\mathcal{U}\rightarrow\mathbb{R}^m,\nonumber\\
\gamma^{\bm{p}}(x,v)&:&\mathbb{R}^m\times\mathcal{U}\rightarrow \mathbb{R}^{m},\forall \bm{p}\in\mathbb{N}^n,\nonumber
\end{eqnarray}
and $\mathcal{U}$ is a nonempty subset of $\mathbb{R}^m$ (control
domain). An admissible control $v(\cdot)$ is a $\mathscr{F}_t-$predictable process with values in $\mathcal{U}$ such that
\begin{eqnarray}
\|v(\cdot)\|=:\sup_{0\leq t\leq T}\left[E|v(t)|^8\right]^{\frac{1}{8}}<\infty
\end{eqnarray}
We denote the set of all admissible controls by $\mathcal{U}_{ad}$. When $\mathcal{U}=\mathbb{R}^m$, we write $L_{\mathscr{F},p}^{\infty,8}[[0,1];\mathbb{R}^m]$ for $\mathcal{U}_{ad}$. The terminal constraint is
\begin{eqnarray}
\mathbb{E}G(x_0,X(T))\in Q\subset\mathbb{R}^k,
\end{eqnarray}
where $G(\cdot,\cdot)=:(G^1(\cdot,\cdot)),\cdots,G^k(\cdot,\cdot)$ and $G^i(\cdot,\cdot):\mathbb{R}^m\times\mathbb{R}^m\rightarrow\mathbb{R}^k$ for $i=1,2,\cdots,k$.

The cost functional is
\begin{eqnarray}
J(v(\cdot),x_0)&=&E\int_0^T\ell(x(t),v(t))dt+Eh(x_0,x(T)),
\end{eqnarray}
where
\begin{eqnarray}
\ell(x,v):\mathbb{R}^m\times\mathcal{U}\rightarrow\mathbb{R},\ \
h(x):\mathbb{R}^m\rightarrow\mathbb{R}.\nonumber
\end{eqnarray}
Our optimal control problem is to find a pair $(y_0,u(\cdot))\in\mathbb{R}^m\times\mathcal{U}_{ad}$ such that (3.1) and (3.3) are satisfied and (3.4) is minimized

Throughout the paper, we make the following assumptions
\begin{assu}
The vector functions $g(x,v)$, $G(y,x)$,$\ell(x,v)$,$h(y,x)$ and $\gamma^{\bm{p}}(x,v)(\bm{p}\in\mathbb{N}_0^n)$ are twice continuously differentiable with respect to $x$, and $G(y,x)$,$h(y,x)$ are differentiable in $y$. They and their derivatives in $x$ or $y$ are continuous in $(x,v)$ and $(y,x)$.The vector functions $g(x,v)$, $G_{y_i}(y,x)$,$G_{x_i}(y,x)$,$\ell_{x_i}(x,v)$,$h_{y_i}(y,x)$,$h_{x_i}(y,x)$,and
$$\left[\sum\limits_{d=1}^{\infty}\sum\limits_{\bm{p}\in\mathbb{N}_d^n}
|\bm{z}^{\bm{p}}(x,v)|^{2k}\right]^{\frac{1}{2k}},\qquad k=1,2,$$
$(i=1,\cdots,n)$,are bounded by $(1+|x|+|y|+|v|)$. The vector functions $G(y,x)$,$\ell(x,v)$,$h(y,x)$ are bounded by $(1+|x|^2+|y|^2+|v|^2)$,$g_{x_i}(x,v)$,$g_{x_ix_j}(x,v)$,$G_{x_ix_j}(y,x)$,$\ell_{x_ix_j}(x,v)$,
$h_{x_ix_j}(y,v)$,and
$$\sum\limits_{d=1}^{\infty}\sum\limits_{\bm{p}\in\mathbb{N}_d^n}
|\bm{z}_{x_i}^{\bm{p}}(x,v)|^{2k},\qquad k=1,2,\qquad\sum\limits_{d=1}^{\infty}\sum\limits_{\bm{p}\in\mathbb{N}_d^n}
|\bm{z}_{x_ix_j}^{\bm{p}}(x,v)|^2$$
$(i,j=1,\cdots,n)$are bounded. Here $x_i,y_i(i=1,\cdots,n)$ stand for the $i$th coordinates of $x$ and $y$ respectively.
\end{assu}

\begin{assu}
The set $Q$ is closed and convex.
\end{assu}

Let
$(y_0,y(\cdot),u(\cdot))$ be an optimal triplet of the problem. For
the given $(x_0,v(\cdot))\in \mathbb{R}^m\times\mathcal{U}_{ad}$, write
$y(\cdot;v(\cdot),x_0)$ for the solution of (3.1). For $v(\cdot)$,
$v_1(\cdot)$, $v_2(\cdot)\in\mathcal{U}_{ad}$, denote
\begin{eqnarray}
\begin{array}{rcl}
\triangle
m(s;v_2,v_1)&\stackrel{\rm{def}}{=}&m(y(s-),v_2)-m(y(s-),v_1),\\
\triangle m(s;v)&\stackrel{\rm{def}}{=}&m(y(s-),v)-m(y(s-),u(s)),\\
m(s;v_1)&\stackrel{\rm{def}}{=}&m(y(s),v_1),\\
m(s)&\stackrel{\rm{def}}{=}&m(y(s),u(s)),
\end{array}
\end{eqnarray}
with $m$ standing for $g$,$\gamma$, $\ell$ and all their
(up to second-) derivatives in $x$.

For $I_0\subset [0,1]$, let $|I_0|$ denote the Lebesgue measure of
the set $I_0$. Let $v(\cdot)$, $v_1(\cdot)$, $v_2(\cdot)\in\mathcal{U}_{ad}$.
Define
\begin{eqnarray}
\hat{d}(v_1(\cdot),v_2(\cdot))\stackrel{\rm{def}}{=}|\{t\in[0,1];E|v_1(\cdot)-v_2(\cdot)|^2>0\}|.
\end{eqnarray}

For $\rho\in (0,T]$, $I_\rho\subset [0,T]$ and $v(\cdot)\in\mathcal{U}_{ad}$,
It is classical to construct a perturbed admissible control in the
following way (spike variation):
\begin{eqnarray}
\begin{array}{lcl}
  u^\rho(s)& \stackrel{\rm{def}}{=} & u(s)\chi_{[0,1]\backslash I_\rho}(s)+v(s)\chi_{I_\rho}(s),\qquad s\in [0,T], \\
  y_0^\rho & \stackrel{\rm{def}}{=} & y_0+|I_\rho|\eta ,\qquad \eta\in \mathbb{R}^n \\
  y^\rho(\cdot)& \stackrel{\rm{def}}{=} & y(\cdot ;u^\rho(\cdot),y_0^\rho),
\end{array}
\end{eqnarray}
with $\chi_A(\cdot)$ denoting the indicator function of some set
$A$. Obviously, we have
\begin{eqnarray}
\hat{d}(u_\rho(\cdot),u(\cdot))=|I_\rho| .
\end{eqnarray}
We can prove that $u^\rho(\cdot)\in\mathcal{U}_{ad}$.

\begin{lem}
Let the Assumption 3.1 hold. Then for $v(\cdot),u(\cdot),u^\rho(\cdot)\in\mathcal{U}_{ad}$
\begin{eqnarray}
\begin{array}{l}
 \sup\limits_{t\in
[0,T]}\mathbb{E}\left|y(t;v(\cdot),x_0)\right|^{8}=O((1+\|v(\cdot)\|)^8),  \\
\sup\limits_{t\in
[0,T]}\mathbb{E}\left|y(t,u(\cdot),y_0)-y(t;u^\rho(\cdot),x_0)\right|^{4}=O(\hat{d}^2(u^\rho(\cdot),u(\cdot))
(1+\|u(\cdot)\|+\|u^\rho(\cdot)\|)^4),  \\
\sup\limits_{t\in
[0,T]}\mathbb{E}\left|y_1(t;u^\rho(\cdot),u(\cdot))\right|^{8}=O(\hat{d}^4(u^\rho(\cdot),u(\cdot))
(1+\|u(\cdot)\|+\|u^\rho(\cdot)\|)^8),  \\
\sup\limits_{t\in
[0,T]}\mathbb{E}\left|y_2(t;u^\rho(\cdot),u(\cdot))\right|^{4}=O(\hat{d}^4(u^\rho(\cdot),u(\cdot))
(1+\|u(\cdot)\|+\|u^\rho(\cdot)\|)^8),  \\
\sup\limits_{t\in [0,T]}\mathbb{E}|y(t;u^\rho,y_0+\hat{d}(u_i,u)\eta)-y(t;u,y_0)-y_1(t;u^\rho,u)-y_2(t;u^\rho,u)|^2\\
=o(\tilde{d}^2(u^\rho(\cdot),u(\cdot))(1+\|u(\cdot)\|+\|u^\rho(\cdot)\|)^8),
\quad as\quad \hat{d}(u^\rho(\cdot),u(\cdot))\rightarrow 0.
\end{array}
\end{eqnarray}
where $y_1(\cdot)$, $y_2(\cdot)$ are the solutions of
\begin{eqnarray}
y_1(t)&=&\int_0^t
g_x(y(s),u(s))y_1(s)ds\nonumber\\
&&+\sum\limits_{d=1}^{\infty}\sum\limits_{\bm{p}\in\mathbb{N}_d^n}\int_0^t\left[
\gamma^{\bm{p}}_x(y(s),u(s))y_1(s)+\triangle\gamma^{\bm{p}}(s,u^\rho(s),u(s))
\right]dH^{\bm{p}}(s)\\
y_2(t)&=&\hat{d}(u^\rho(\cdot),u(\cdot))\eta\nonumber\\
&&+\int_0^t\left[g_x(y(s),u(s))y_2(s)+\triangle
g(s,u^\rho(s),u(s))+\frac{1}{2}g_{xx}(y(s),u(s))y_1(s)y_1(s)\right]ds\nonumber\\
&&+\sum\limits_{d=1}^{\infty}\sum\limits_{\bm{p}\in\mathbb{N}_d^n}\int_0^t\left[\gamma^{\bm{p}}_x(y(s),u(s))y_2(s)
+\frac{1}{2}\gamma^{\bm{p}}_{xx}(y(s),u(s))y_1(s)y_1(s)\right]dH^{\bm{p}}(s)\nonumber\\
&&+\sum\limits_{d=1}^{\infty}\sum\limits_{\bm{p}\in\mathbb{N}_d^n}\int_0^t\triangle\gamma^{\bm{p}}_x(s,u^\rho(s),u(s))
y_1(x)dH^{\bm{p}}(s)
\end{eqnarray}
where $f_{xx}yy=\sum_{i,j=1}^mf_{x^ix^j}y^iy^j$ for $f=g,\gamma^{\bm{p}}$.
\end{lem}

Proof. Without loss of generality, we assume $\eta=0$. Define
\begin{eqnarray}
\int_{I_\rho}g_0(s)dH^{\bm{p}}(s)&=:&\int\chi_{I_\rho}(s)g_0(s)dH^{\bm{p}}(s).\quad\forall\bm{p}\in\mathbb{N}^n \nonumber
\end{eqnarray}
We have the following inequalities for $p>1$:
\begin{eqnarray}
\begin{array}{rcl}
\mathbb{E}\left|\int_{I_\rho}f_0(s)ds\right|^p&\leq&C_p|I_\rho|^{p-1}\mathbb{E}\int_{I_\rho}|f_0(s)|^pds,\\
\mathbb{E}\left|\int_{I_\rho}g_0(s)dH^{\bm{p}}(s)\right|^{2p}&\leq&C_p|I_\rho|^{p-1}\mathbb{E}
\int_{I_\rho}\left|g_0(s,z)\right|^{2p}ds,\quad\forall\bm{p}\in\mathbb{N}^n.
\end{array}
\end{eqnarray}

By virtue of the Assumption 3.1, we have
\begin{eqnarray}
\begin{array}{l}
 \sup\limits_{t\in
[0,T]}\mathbb{E}\left|y(t)\right|^{8}=O((1+\|v(\cdot)\|+\|u(\cdot)\|)^8),  \\
\sup\limits_{t\in
[0,T]}\mathbb{E}\left|\triangle g(t,u^\rho(s)))\right|^{4}=O((1+\|u^\rho(\cdot)\|+\|u(\cdot)\|)^4),  \\
\sup\limits_{t\in
[0,T]}\mathbb{E}\left|\triangle \gamma^{\bm{p}}(t;u^\rho(s))\right|^8=O((1+\|u^\rho(\cdot)\|+\|u(\cdot)\|)^8),\quad\forall\bm{p}\in\mathbb{N}^n.
\end{array}
\end{eqnarray}

Then we can obtain the following inequalities by using (3.12):
\begin{eqnarray}
\begin{array}{l}
\mathbb{E}\left|\int_0^T\triangle g(t,u^\rho(s)))\right|^{4}=O(|I_\rho|^4(1+\|v(\cdot)\|+\|u(\cdot)\|)^4),  \\
\mathbb{E}\left|\int_0^T\triangle \gamma^{\bm{p}}(t;u^\rho(\cdot))\right|^{8}=O(|I_\rho|^4(1+\|v(\cdot)\|+\|u(\cdot)\|)^8),
\quad\forall\bm{p}\in\mathbb{N}^n.
\end{array}
\end{eqnarray}
Then the first four estimates of (3.9) are easily proved by using the familiar elementary inequalities
\begin{eqnarray}
(m_1+m_2)^i&\leq&C(|m_1|^i+|m_2|^i),i=4,8\nonumber
\end{eqnarray}
and the well-known Gronwall's inequality.

The proof for the last estimate follows. Set $ y_3=y_1+y_2$. We have
\begin{eqnarray}
&&\int_0^t
g(y+y_3,u^{\rho})ds+\sum\limits_{d=1}^{\infty}\sum\limits_{\bm{p}\in\mathbb{N}_d^n}
\int_0^t\gamma^{\bm{p}}(y+y_3,u^{\rho})dH^{\bm{p}}(s),\nonumber\\
&=&\int_0^t\left[g(y,u^{\rho})+g_x(y,u^\rho)y_3+\int_0^1\int_0^1\lambda
g_{xx}(y+\lambda\mu y_3,u^\rho)d\lambda d\mu y_3y_3\right]ds\nonumber\\
&&+\sum\limits_{d=1}^{\infty}\sum\limits_{\bm{p}\in\mathbb{N}_d^n}
\int_0^t\left[\gamma^{\bm{p}}(y,u^{\rho})+\gamma^{\bm{p}}_x(y,u^\rho)y_3+\int_0^1\int_0^1\lambda
\gamma^{\bm{p}}_{xx}(y+\lambda\mu y_3,u^\rho)d\lambda d\mu y_3y_3\right]dH^{\bm{p}}(s)\nonumber\\
&=&\int_0^t g(y,u)ds+\int_0^t
g_x(y,u)y_3ds+\int_0^t\triangle g(s,u^\rho(s),u(s))ds\nonumber\\
&&+\int_0^t\triangle g_x(s,u^\rho(s),u(s))y_3(s)ds+\int_0^t\frac{1}{2}g_{xx}(y,u)y_3(s)y_3(s)ds\nonumber\\
&&+\int_0^t\int_0^1\int_0^1\lambda\left[g_{xx}(y+\lambda\mu
y_3,u^\rho)-g_{xx}(y,u)\right]d\lambda d\mu
y_3y_3ds\nonumber\\
&&+\sum_{d=1}^{\infty}\sum\limits_{\bm{p}\in\mathbb{N}_d^n}\int_0^t \gamma^{\bm{p}}(y,u)dH^{\bm{p}}(s)+
\sum\limits_{d=1}^{\infty}\sum\limits_{\bm{p}\in\mathbb{N}_d^n}\int_0^t \gamma^{\bm{p}}_x(y,u)y_3 dH^{\bm{p}}(s)\nonumber\\
&&+\sum_{d=1}^{\infty}\sum\limits_{\bm{p}\in\mathbb{N}_d^n}\int_0^t
\gamma^{\bm{p}}(s,u^\rho(s),u(s))dH^{\bm{p}}(s)+
\sum\limits_{d=1}^{\infty}\sum\limits_{\bm{p}\in\mathbb{N}_d^n}
\int_0^t\triangle\gamma^{\bm{p}}_x(s,u^\rho(s),u(s))y_3dH^{\bm{p}}(s)\nonumber\\
&&+\sum\limits_{d=1}^{\infty}\sum\limits_{\bm{p}\in\mathbb{N}_d^n}\int_0^t
\frac{1}{2}\gamma^{\bm{p}}_{xx}(y,u)y_3(s)y_3(s)dH^{\bm{p}}(s)\nonumber\\
&&+\sum\limits_{d=1}^{\infty}\sum\limits_{\bm{p}\in\mathbb{N}_d^n}\int_0^t\int_0^1\int_0^1\lambda
\left[\gamma^{\bm{p}}_{xx}(y+\lambda\mu y_3,u^\rho)-\gamma^{\bm{p}}_{xx}(y,u)\right]y_3y_3d\lambda
d\mu dH^{\bm{p}}(s)\nonumber\\
&=&y(t)+y_3(t)-y_0+\int_0^tG^{\rho}(s)ds+\sum\limits_{d=1}^{\infty}\sum\limits_{\bm{p}\in\mathbb{N}_d^n}\int_0^t \Xi^{\rho,\bm{p}}(s)
dH^{\bm{p}}(s),\nonumber
\end{eqnarray}
where
\begin{eqnarray}
G^{\rho}(s)&=&\frac{1}{2}g_{xx}(y(s),u(s))(y_2(s)y_2(s)+2y_1(s)y_2(s))\nonumber\\
&&+\triangle g_x(y(s),u^\rho(s),u(s))y_2(s)\nonumber\\
&&+\int_0^1\int_0^1\lambda\left[g_{xx}(y+\lambda\mu
y_3,u^\rho)-g_{xx}(y,u)\right]d\lambda d\mu
y_3(s)y_3(s)\nonumber\\
\Xi^{\rho,\bm{p}}(s)&=&\frac{1}{2}\gamma^{\bm{p}}_{xx}(y(s),u(s))(y_2(s)y_2(s)+2y_1(s)y_2(s))\nonumber\\
&&+\triangle\gamma^{\bm{p}}_x(y(s),u^\rho(s),u(s))y_2(s)\nonumber\\
&&+\int_0^1\int_0^1\lambda\left[\gamma^{\bm{p}}_{xx}(y+\lambda\mu
y_3,u^\rho)-\gamma^{\bm{p}}_{xx}(y,u)\right]d\lambda d\mu
y_3(s)y_3(s)\nonumber
\end{eqnarray}
Since
\begin{eqnarray}
y(t)+y_3(t)&=&y_0+\int_0^t
g(y+y_3,u^{\rho})ds+\sum_{d=1}^{\infty}\sum\limits_{\bm{p}\in\mathbb{N}_d^n}
\int_0^t\gamma^{\bm{p}}(y+y_3,u^{\rho})dH^{\bm{p}}(s)\nonumber\\
&&-\int_0^tG^{\rho}(s)ds-\sum\limits_{d=1}^{\infty}\sum\limits_{\bm{p}\in\mathbb{N}_d^n}\int_0^t\Xi^{\rho,\bm{p}}(s)dH^{\bm{p}}(s)\nonumber
.
\end{eqnarray}
and
\begin{eqnarray}
y^\rho(t)=y_0+\int_0^t
g(y^\rho(s),u^{\rho}(s))ds+\sum_{d=1}^{\infty}\sum\limits_{\bm{p}\in\mathbb{N}_d^n}\int_0^t\gamma^{\bm{p}}(y^\rho(s),u^{\rho}(s))dH^{\bm{p}}(s),\nonumber
\end{eqnarray}
we can derive that
\begin{eqnarray}
(y^\rho-y-y_3)(t)&=&\int_0^tA^\rho(s)(y^\rho-y-y_3)(s)ds\nonumber\\
&&+\sum_{d=1}^{\infty}\sum\limits_{\bm{p}\in\mathbb{N}_d^n}\int_0^tF^{\rho,\bm{p}}(s)(y^\rho-y-y_3)(s)dH^{\bm{p}}(s)\nonumber\\
&&+\int_0^tG^{\rho}(s)ds+\sum\limits_{d=1}^{\infty}\sum\limits_{\bm{p}\in\mathbb{N}_d^n}\int_0^t
\Xi^{\rho,\bm{p}}(s)dH^{\bm{p}}(s).\nonumber
\end{eqnarray}
\begin{eqnarray}
\left|A^\rho(s,\omega)\right|+\sum_{d=1}^{\infty}\sum\limits_{\bm{p}\in\mathbb{N}_d^n}\left|F^{\rho,\bm{p}}(s,\omega)\right|\leq C \quad\forall s,\ \ \forall \omega . \nonumber
\end{eqnarray}
and
\begin{eqnarray}
\sup_{0\leq t\leq
T}E\left(\left|\int_0^tG^\rho(s)ds\right|^2+\sum_{d=1}^{\infty}\sum\limits_{\bm{p}\in\mathbb{N}_d^n}\left|\int_0^t
\Xi^{\rho,\bm{p}}(s)dH^{\bm{p}}(s)\right|^2\right)=o(|I_\rho|^2(1+\|u^\rho(\cdot)\|+\|u(\cdot)\|)^8).\nonumber
\end{eqnarray}
From these we can use It$\hat{o}$'s formula and
Gronwall's inequality to obtain the fifth estimate (3.9). The proof is
completed.$\Box$

\begin{lem}
Assume that $l(\cdot)$ is a scalar-valued Lebesgue integrable
function defined on $[0,T]$. Then for $\rho\in (0,T]$, there exists
a measurable subset $I_\rho\subset [0,T]$, such that
\begin{eqnarray}
\begin{array}{rcl}
|I_\rho|&=&\rho,\\
\int_{I_\rho}l(s)ds&=&\rho\int_{[0,T]}l(s)ds+o(\rho),  \qquad
\rho\rightarrow 0 .
\end{array}
\end{eqnarray}
\end{lem}

The proof is quite elementary and the reader is referred to [18].

\section{Adjoint equations and the maximum principle}

The Hamiltonian is defined as
\begin{eqnarray}
H(x,v,\lambda,p,J)=\lambda
\ell(x,v)+(p,g(x,v))+\sum_{i=1}^{\infty}\sum\limits_{\bm{p}\in\mathbb{N}_i^r}
(J^{\bm{p}},\gamma^{\bm{p}}(x,v))\
. \nonumber
\end{eqnarray}
this is a map from $\mathbb{R}^m\times\mathcal{U}\times\mathbb{R}\times\mathbb{R}^m\times (M_T^2(l^2))^m$ into $\mathbb{R}$. Here we have used $(\cdot,\cdot)$ for the scalar product of Euclidean spaces.

From Lemma 2.2 and Assumption 3.1, we see for the given $p(T)\in L^2(\Omega,\mathscr{F}_T;\mathbb{R}^m)$, $P(T)\in L^2(\Omega,\mathscr{F}_T;\mathbb{R}^{m\times m})$ that the It$\hat{o}$-type adjoint equations
\begin{eqnarray}
-dp(t)&=&\left[g_x^\top(y(t),u(t))p(t)+\sum_{d=1}^{\infty}\sum\limits_{\bm{p}\in\mathbb{N}_d^n}\gamma_{x}^{\bm{p}}(y(t),u(t))^\top
J^{\bm{p}}(t)+\lambda\ell_x(y(t),u(t))\right]dt\nonumber\\
 &&-\sum_{d=1}^{\infty}\sum\limits_{\bm{p}\in\mathbb{N}_d^n}J^{\bm{p}}(t)dH^{\bm{p}}(t)\nonumber\\
 p(T)&=&h_x(y(T)).
\end{eqnarray}
and
\begin{eqnarray}
-dP(t)&=&\left[g_x^\top(y(t),u(t))P(t)+P(t)g_x(y(t),u(t))+\sum_{d=1}^{\infty}\sum\limits_{\bm{p}\in\mathbb{N}_d^n}\gamma_{x}^{\bm{p}}(y(t),u(t))^\top
P(t)\gamma_{x}^{\bm{p}}(y(t),u(t))\right.\nonumber\\
&&+\sum_{d=1}^{\infty}\sum\limits_{\bm{p}\in\mathbb{N}_d^n}\gamma_{x}^{\bm{p}}(y(t),u(t))R^{\bm{p}}(t)\nonumber\\
&&\left.+\sum_{d=1}^{\infty}\sum\limits_{\bm{p}\in\mathbb{N}_d^n}^\top R^{\bm{p}}(t)\gamma_{x}^{\bm{p}}(y(t),u(t))+H_{xx}(y(t),u(t),\lambda,p(t),J(t))\right]dt\nonumber\\
&&-\sum_{d=1}^{\infty}\sum\limits_{\bm{p}\in\mathbb{N}_d^n}R^{\bm{p}}(t)dH^{\bm{p}}(t)\nonumber\\
P(T)&=&h_{xx}(y(T))
\end{eqnarray}
admit unique solutions $(p(\cdot),\{J^{\bm{p}}(\cdot)\}_{\bm{p}\in\mathbb{N}^n})$ and
$(P(\cdot),\{R^{\bm{p}}(\cdot)\}_{\bm{p}\in\mathbb{N}^n})$, with $p(\cdot)$ and $P(\cdot)$ being cadlag processes.

Define the following function:
\begin{eqnarray}
\Phi(s,z;\varepsilon)\stackrel{\rm{def}}{=}\inf_{(t,\bar{z})\in
(-\infty,J(u(\cdot),y_0)-\varepsilon]\times
Q}\sqrt{|t-s|^2+|\bar{z}-z|^2}
\end{eqnarray}
Tang and Li [34] had proved the following result.
\begin{lem}
For given $\varepsilon>0$, the function $\Phi(s,z;\varepsilon)$ is
continuously differentiable on the open set
$\hat{Q}\stackrel{\rm{def}}{=}\{(s,z):\Phi(s,z;\varepsilon)>0\}$.
Moreover, when $\Phi(s,z;\varepsilon)>0$, we have
\begin{eqnarray}
\begin{array}{rcl}
<\Phi_z(s,z;\varepsilon),\hat{z}-z>&\leq& 0, \forall
\hat{z}\in Q,\\
\Phi_s(s,z;\varepsilon)&\geq&0,\\
|\Phi_s(s,z;\varepsilon)|^2+|\Phi_z(s,z;\varepsilon)|^2&=&1 .
\end{array}
\end{eqnarray}
\end{lem}

They introduce the smooth function $\alpha(\cdot)$ defined by
\begin{eqnarray}
\alpha(t,z)\stackrel{\rm{def}}{=}\left\{\begin{array}{cc}
                                        C exp(t^2+|z|^2-1)^{-1},& t^2+|z|^2<1, \\
                                        0, &t^2+|z|^2\geq 1.
                                      \end{array}
\right.\nonumber
\end{eqnarray}
Choose the constant $C$ such that
\begin{eqnarray}
\int_{\mathbb{R}\times\mathbb{R}^k}\alpha(t,z)dtdz=1.\nonumber
\end{eqnarray}
Set
\begin{eqnarray}
\alpha_\delta(t,z)=\delta^{-(k+1)}\alpha\left(\frac{t}{\delta},\frac{z}{\delta}\right).
\end{eqnarray}
They also define the smooth approximation
$\Psi(\cdot,\cdot;\varepsilon,\delta)$ of
$\Phi(\cdot,\cdot;\varepsilon)$ as follows:
\begin{eqnarray}
\Psi(s,z;\varepsilon,\delta)\stackrel{\rm{def}}{=}\int_{\mathbb{R}\times\mathbb{R}^k}
\Phi(s-\bar{s},z-\bar{z};\varepsilon)\alpha_\delta(\bar{s},\bar{z})d\bar{s}d\bar{z}=1.
\end{eqnarray}

Then it is easy to have
\begin{eqnarray}
0\leq\Psi(J(u(\cdot),y_0),EG(y_0,y(T));\varepsilon,\delta)\leq\varepsilon+\sqrt{2}\delta\nonumber
\end{eqnarray}
Moreover, Tang and Li [34] gave the following lemma.

\begin{lem}
For $\hat{Q}$ defined in Lemma 4.1, we have for $(s,z)\in \hat{Q}$,
\begin{eqnarray}
\lim_{\delta\rightarrow
0+}\Psi_s(s,z;\varepsilon,\delta)&=&\Phi_s(s,z;\varepsilon),\nonumber\\
\lim_{\delta\rightarrow
0+}\Psi_z(s,z;\varepsilon,\delta)&=&\Phi_z(s,z;\varepsilon) .
\end{eqnarray}
\end{lem}
Our main result in this paper is almost similar to that in Tang and Li [34] in many places:
\begin{theorem}
Assume Assumptions 3.1 and 3.2 hold. Let $(y_0,y(\cdot),u(\cdot))$ be an
optimal triplet. Then there exist
\begin{eqnarray}
\begin{array}{l}
 0\leq \lambda\in \mathbb{R}, \quad
\mu\stackrel{\rm{def}}{=}\{\mu^i\}_1^k\in\mathbb{R}^k,  \\
\left(p(\cdot),J(\cdot)\right)\in
L_{\mathcal{F}}^2(0,T;\mathbb{R}^m)\times
L_{\mathcal{F}}^2(0,T;(M_T^2(l^2))^m)   \\
 \left(P(\cdot),R(\cdot)\right)\in
L_{\mathcal{F}}^2(0,T;\mathbb{R}^{m\times m})\times
L_{\mathcal{F}}^2(0,T;(M_T^2(l^2))^{m\times m})
\end{array}\nonumber
\end{eqnarray}
such that we have the following.

1) The nontrivial condition
\begin{eqnarray}
|\lambda|^2+|\mu|^2=1,
\end{eqnarray}
is satisfied.

2) The It$\hat{o}$-type adjoint equations (4.1),(4.2), as well as
\begin{eqnarray}
\left\{\begin{array}{ccc}
         p(T)&=&\lambda h_x(y_0,y(T))+\sum\limits_{j=1}^k\mu^jG_x^j(y_0,y(T)), \\
         p(0)&=&-\lambda Eh_x(y_0,y(T))-\sum\limits_{j=1}^k\mu^jEG_y^j(y_0,y(T))
       \end{array}\right.
\end{eqnarray}
and
\begin{eqnarray}
p(T)=\lambda h_{xx}(y_0,y(T))+\sum_{j=1}^k\mu^jG_{xx}^j(y_0,y(T)),
\end{eqnarray}
are satisfied, with $p(\cdot)$ and $P(\cdot)$ being cadlag
processes.

3) The following maximum condition holds:
\begin{eqnarray}
\begin{array}{l}
H(y(s-),v,\lambda,p(s-),J(s))-H(y(s-),u(s),\lambda,K(s),J(s))\\
+\frac{1}{2}tr
P(s-)\left[\sum\limits_{d=1}^\infty\sum\limits_{\bm{p}\in\mathbb{N}_d^n}
\triangle\gamma^{\bm{p}}(s;v)\sum\limits_{d=1}^\infty\sum\limits_{\bm{p}\in\mathbb{N}_d^n}
\triangle\gamma^{\bm{p}\top}(s;v)\right]\\
+\frac{1}{2}tr\left[ \sum\limits_{d=1}^\infty\sum\limits_{\bm{p}\in\mathbb{N}_d^n}
R^{\bm{p}}(s)\right]^\top\left[\sum\limits_{d=1}^\infty\sum\limits_{\bm{p}\in\mathbb{N}_d^n}
\triangle\gamma^{\bm{p}}(s;v)\sum\limits_{d=1}^\infty\sum\limits_{\bm{p}\in\mathbb{N}_d^n}
\triangle\gamma^{\bm{p}\top}(s;v)\right]\\
\geq 0,\quad\forall v(\cdot)\in \mathcal{U},\quad a.e.a.s.;
\end{array}
\end{eqnarray}

4) The following transversality condition holds:
\begin{eqnarray}
<\mu,z-EG(y_0,y(T))>\leq 0,\quad \forall z\in Q .
\end{eqnarray}
\end{theorem}

\noindent \textbf{Proof}

\textsl{Step 1. Applying Ekeland's variational principle.} We first
consider the case that the set $\mathcal{U}_{ad}$ is bounded in $L_{\mathscr{F},p}^{\infty,8}[[0,T];\mathbb{R}^m]$; the
unbounded case can be reduced to the bounded case. Assume that
\begin{eqnarray}
\mathcal{U}_{ad} \quad is \quad bounded \quad in \quad L_{\mathscr{F},p}^{\infty,8}[[0,T];\mathbb{R}^m]
\end{eqnarray}

An applications of Ekeland's variational principle will lead to the
reduction of a general end-constraint problem to a family of free
end-constraint problems.

Define the following auxiliary function
\begin{eqnarray}
J(v(\cdot),x_0;\varepsilon,\delta)=\Psi(J(v(\cdot),x_0),EG(x_0,x(T));\varepsilon,\delta)
\end{eqnarray}
with $\Psi(\cdot,\cdot;\varepsilon,\delta)$ being defined as in (4.6). Then consider the metric space $(\mathbb{R}^m\times\mathcal{U}_{ad},d)$
with the distance $d$ defined by
\begin{eqnarray}
d((x_1,v_1(\cdot)),(x_2,v_2(\cdot)))=\sqrt{|x_1-x_2|^2+\hat{d}^2(v_1(\cdot),v_2(\cdot))}.
\end{eqnarray}
Tang and Li [34] verify that $\Psi(\cdot,\cdot;\varepsilon,\delta)$ is
complete and $J(v(\cdot),x_0;\varepsilon,\delta)$ is continuous and
bounded. Also, we have for any given $\varepsilon>0$,
\begin{eqnarray}
\begin{array}{l}
  \Phi(J(v(\cdot),x_0),\quad EG(x_0,x(T));\varepsilon)>0, \quad \forall
(x_0,v(\cdot))\in \mathbb{R}^m\times \mathcal{U}_{ad}; \\
 \Phi(J(v(\cdot),y_0),\quad EG(y_0,y(T));\varepsilon)=\varepsilon; \\
  J(v(\cdot),x_0;\varepsilon,\delta)>0, \quad \forall
(x_0,v(\cdot))\in \mathbb{R}^m\times \mathcal{U}_{ad}, \\
 \quad for \quad
sufficiently \quad small \quad \delta>0;  \\
J(u(\cdot),y_0;\varepsilon,\delta)\leq
\varepsilon+2\delta+\inf_{(x_0,v(\cdot))\in\mathbb{R}^n\times
\mathcal{U}_{ad}}J(v(\cdot),x_0;\varepsilon,\delta)
\end{array}
\end{eqnarray}
Therefore we can apply Ekeland's variational principle (cf.[10]) and
conclude that there exist $u^{\varepsilon\delta}\in\mathcal{U}_{ad}$ and
$y_0^{\varepsilon\delta}\in \mathbb{R}^m$ such that
\begin{eqnarray}
\begin{array}{ll}
1)&J(u^{\varepsilon\delta}(\cdot),y_0^{\varepsilon\delta};\varepsilon,\delta)\leq
\varepsilon+2\delta;\\
2)&d((y_0^{\varepsilon\delta},u^{\varepsilon\delta}(\cdot)),(y_0,u(\cdot)))\leq\sqrt{\varepsilon+2\delta}\\
3)&\bar{J}(v(\cdot),x_0;\varepsilon,\delta)\stackrel{\rm{def}}{=}J(v(\cdot),x_0;\varepsilon,\delta)
+\sqrt{\varepsilon+2\delta}d((x_0,v(\cdot)),(y_0^{\varepsilon\delta},u^{\varepsilon\delta}(\cdot)))\\
&\geq J(u^{\varepsilon\delta}(\cdot),y_0^{\varepsilon\delta}),\quad
\forall (x_0,v(\cdot))\in \mathbb{R}^m\times\mathcal{U}_{ad} .
\end{array}
\end{eqnarray}

Set
\begin{eqnarray}
\begin{array}{ll}
\lambda^{\epsilon\delta}&\stackrel{\rm{def}}{=}\Psi_s(J(u^{\epsilon\delta}(\cdot),y_0^{\epsilon\delta}),EG(y_0^{\epsilon\delta},
y^{\epsilon\delta}(T));\varepsilon,\delta),\\
\mu^{\epsilon\delta}&\stackrel{\rm{def}}{=}\Psi_z(J(u^{\epsilon\delta}(\cdot),y_0^{\epsilon\delta}),EG(y_0^{\epsilon\delta},
y^{\epsilon\delta}(T));\varepsilon,\delta).
\end{array}
\end{eqnarray}

and
\begin{eqnarray}
\begin{array}{ll}
\lambda^\varepsilon\stackrel{\rm{def}}{=}\lambda^{\varepsilon\delta(\varepsilon)},&
\mu^\varepsilon\stackrel{\rm{def}}{=}\mu^{\varepsilon\delta(\varepsilon)},
\\
y_0^\varepsilon\stackrel{\rm{def}}{=}y_0^{\varepsilon\delta(\varepsilon)},&
u^\varepsilon(\cdot)\stackrel{\rm{def}}{=}u^{\varepsilon\delta(\varepsilon)}(\cdot).
\end{array}\nonumber
\end{eqnarray}
Tang and Li [34] showed that for each sufficiently small $\varepsilon>0$, we can choose $\delta(\varepsilon)>0$ such that $\lambda^\varepsilon\geq 0$ and
$\mu^\varepsilon\in \mathbb{R}^k$ satisfy the following:
\begin{eqnarray}
\begin{array}{rl}
\lim\limits_{\delta\rightarrow
0+}(|\lambda^{\epsilon}|^2+|\mu^{\epsilon}|^2)&=1,\\
<\mu^\varepsilon,z-EG(y_0^\varepsilon,y^\varepsilon(T))>&\leq
\delta(\varepsilon)\leq \varepsilon .
\end{array}
\end{eqnarray}

\textsl{Step 2. Computing the first-order component of the cost
variation.} In this and the next steps, we look for the necessary
conditions for the minimization of
$\bar{J}(v(\cdot),x_0;\varepsilon,\delta)$ at
$(y_0^\varepsilon,u^\varepsilon(\cdot))$.

For given $(\eta,v(\cdot))\in \mathbb{R}^m\times\mathcal{U}_{ad}$, set
\begin{eqnarray}
\begin{array}{rcl}
u^{\varepsilon\rho}(t)&=&u^\varepsilon(t)\chi_{[0,1]\backslash I_\rho}(t)+v(t)\chi_{I_\rho}(t),\\
y_0^{\varepsilon\rho}&=&y_0^\varepsilon+|I_\rho|\eta,\\
y^{\varepsilon\rho}(\cdot)&=&y(\cdot;u^{\varepsilon\rho}(\cdot),y_0^{\varepsilon\rho}).
\end{array}
\end{eqnarray}

We introduce, as in (3.4), the following simplified notations:
\begin{eqnarray}
\begin{array}{rcl}
\triangle m^\varepsilon(s;v)&\stackrel{\rm{def}}{=}&
m(y^\varepsilon(s),v)-m(y^\varepsilon(s),u^\varepsilon(s)),\\
m^\varepsilon(s)&\stackrel{\rm{def}}{=}&
m(y^\varepsilon(s),u^\varepsilon(s)),
\end{array}
\end{eqnarray}
with $m$ standing for $g$,$\gamma$, $\ell$ and all their
(up to second-) derivatives in $x$.

Let $y^{\varepsilon\rho}(\cdot)$ be the solution of (3.1) corresponding
to $(y_0^{\varepsilon\rho},u^{\varepsilon\rho}(\cdot))$. We define,
as in (3.9) and (3.10), the half- and first-order processes
$y_1^\varepsilon(\cdot)$, $y_2^\varepsilon(\cdot)$, respectively, by
\begin{eqnarray}
y_1^\varepsilon(t)&=&\int_0^tg_x(y^\varepsilon(s),u^\varepsilon(s))y_1^\varepsilon(s)
ds\nonumber\\
&&+\sum_{d=1}^{\infty}\sum\limits_{\bm{p}\in\mathbb{N}_d^n}\int_0^t\left[
\gamma^{\bm{p}}_x(y^\varepsilon(s),u^\varepsilon(s))y_1^\varepsilon(s)+\triangle\gamma^{\varepsilon,\bm{p}}
(s;u^{\varepsilon\rho}(s))\right]dH^{\bm{p}}(s)
\end{eqnarray}
and
\begin{eqnarray}
&&y_2^\varepsilon(t)\nonumber\\
&=&\int_0^t\left[g_x(y^\varepsilon(s),u^\varepsilon(s))y_2^\varepsilon(s)
+\triangle g^\varepsilon(s;u^{\varepsilon\rho}(s))
+\frac{1}{2}g_{xx}(y^\varepsilon(s),u^\varepsilon(s))y_1^\varepsilon(s)y_1^\varepsilon(s)\right]ds\nonumber\\
&&+\sum_{d=1}^{\infty}\sum\limits_{\bm{p}\in\mathbb{N}_d^n}\int_0^t\left[\gamma^{\bm{p}}_x(y^\varepsilon(s),u^\varepsilon(s))y_2^\varepsilon(s)
+\frac{1}{2}\gamma^{\bm{p}}_{xx}(y^\varepsilon(s),u^\varepsilon(s))y_1^\varepsilon(s)y_1^\varepsilon(s)\right]dH^{\bm{p}}(s)\nonumber\\
&&+\int_0^t\triangle
g_x^\varepsilon(s;u^{\varepsilon\rho}(s))y_1^\varepsilon(s)ds
+\sum_{d=1}^{\infty}\sum\limits_{\bm{p}\in\mathbb{N}_d^n}\int_0^t\triangle\gamma^{\varepsilon,\bm{p}}_x(s,u^{\varepsilon\rho}(s))
y_1^\varepsilon(x)dH^{\bm{p}}(s)+|I_\rho|\eta
\end{eqnarray}
From Lemma 3.1, we can have
\begin{eqnarray}
\begin{array}{l}
\sup\limits_{0\leq t\leq T}E|y_1^\varepsilon(t)|^8=O(|I_\rho|^4),\\
\sup\limits_{0\leq t\leq T}E|y_2^\varepsilon(t)|^8=O(|I_\rho|^4),\\
\sup\limits_{0\leq t\leq
T}E|y^{\varepsilon\rho}(t)-y^\varepsilon(t)-y_1^\varepsilon(t)-y_2^\varepsilon(t)|^2=o(|I_\rho|^4),\\
\qquad \qquad \qquad \qquad \qquad \qquad as \quad
|I_\rho|\rightarrow 0.
\end{array}
\end{eqnarray}
In this step, we are to calculate the first-order component of the
cost variation.

From 3) in (4.17) of Step 1, we have
\begin{eqnarray}
\begin{array}{cl}
&-|I_\rho|\sqrt{\varepsilon+2\delta}\sqrt{1+|\eta|^2}\\
\leq&J(u^{\epsilon\rho}(\cdot),y^{\epsilon\rho}(0);\varepsilon)-J(u^{\epsilon}(\cdot),y_0^{\epsilon};\varepsilon)\\
\leq&\lambda^\varepsilon
[J(u^{\epsilon\rho}(\cdot),y_0^{\epsilon}+|I_\rho|\eta)-J(u^\varepsilon(\cdot),y_0^\varepsilon)]\\
&+\sum\limits_{j=1}^m\mu^{\varepsilon
j}[EG^j(y_0^\varepsilon+|I_\rho|\eta,y^{\varepsilon\rho}(T))-EG^j(y_0^\varepsilon,y^\varepsilon(T))]\\
&+O(|J(u^{\epsilon\rho}(\cdot),y_0^{\epsilon}+|I_\rho|\eta)-J(u^\varepsilon(\cdot),y_0^\varepsilon)|^2)\\
&+\sum\limits_{j=1}^mO(|EG^j(y_0^\varepsilon+|I_\rho|\eta,y^{\varepsilon\rho}(T))-EG^j(y_0^\varepsilon,y^\varepsilon(T))|^2)
\end{array}
\end{eqnarray}

Using (4.24), we have
\begin{eqnarray}
\begin{array}{cl}
&J(u^{\epsilon\rho}(\cdot),y_0^{\epsilon}+|I_\rho|\eta)-J(u^\varepsilon(\cdot),y_0^\varepsilon)\\
=&|I_\rho|<Eh_y(y_0^\varepsilon,y^\varepsilon(T)),\eta>
+E<h_x(y_0^\varepsilon,y^\varepsilon(T)),y_1^\varepsilon(T)+y_2^\varepsilon(T)>\\
&+\frac{1}{2}Ey_1^{\varepsilon
\top}(T)h_{xx}(y_0^\varepsilon,y^\varepsilon(T))y_1^\varepsilon(T)\\
&+E\int_0^T\ell_x(y^\varepsilon(s),u^\varepsilon(s))[y_1^\varepsilon(s)+y_2^\varepsilon(s)]ds
+\frac{1}{2}E\int_0^Ty_1^{\varepsilon
\top}(s)\ell_{xx}(y^\varepsilon(s),u^\varepsilon(s))y_1^\varepsilon(s)ds\\
&+E\int_0^T\triangle\ell^\varepsilon(s,u^{\varepsilon\rho}(s))ds
+o(|I_\rho|)
\end{array}
\end{eqnarray}
and similarly
\begin{eqnarray}
\begin{array}{cl}
&EG^j(y_0^\varepsilon+|I_\rho|\eta,y^{\varepsilon\rho}(T))-EG^j(y_0^\varepsilon,y^\varepsilon(T))\\
=&|I_\rho|<EG^j_y(y_0^\varepsilon,y^\varepsilon(T)),\eta>_n
+E<G_x^j(y_0^\varepsilon,y^\varepsilon(T)),y_1^\varepsilon(T)+y_2^\varepsilon(T)>\\
&+\frac{1}{2}Ey_1^{\varepsilon
\top}(T)G_{xx}^j(y_0^\varepsilon,y^\varepsilon(T))y_1^\varepsilon(T)+o(|I_\rho|)
\end{array}
\end{eqnarray}

From Lemma 2.2, we see that
\begin{eqnarray}
-dp^\varepsilon(t)&=&\left[g_x^\top(y^\varepsilon(t),u^\varepsilon(t))p^\varepsilon(t)\right.\nonumber\\
&&\left.+\sum_{d=1}^{\infty}\sum\limits_{\bm{p}\in\mathbb{N}_d^n}\gamma_{x}^{\bm{p}}(y^\varepsilon(t),u^\varepsilon(t))^\top
J^{\bm{p},\varepsilon}(t)
 +\lambda^\varepsilon\ell_x(y^\varepsilon(t),u^\varepsilon(t))\right]dt\nonumber\\
 &&-\sum_{d=1}^{\infty}\sum\limits_{\bm{p}\in\mathbb{N}_d^n}J^{\bm{p},\varepsilon}(t)dH^{\bm{p}}(t)\nonumber\\
 p^\varepsilon(T)&=&\lambda^\varepsilon h_x(y_0^\varepsilon,y^\varepsilon(T))
 +\sum_{j=1}^k\mu^{\varepsilon j}G_x^j(y_0^\varepsilon,y^\varepsilon(T)).
\end{eqnarray}
and
\begin{eqnarray}
-dP^\varepsilon(t)&=&\left[g_x^\top(y^\varepsilon(t),u^\varepsilon(t))P^\varepsilon(t)
+P^\varepsilon(t)g_x(y^\varepsilon(t),u^\varepsilon(t))\right.\nonumber\\
&&+\sum_{d=1}^{\infty}\sum\limits_{\bm{p}\in\mathbb{N}_d^n}\gamma_{x}^{\bm{p}}(y^\varepsilon(t),u^\varepsilon(t))^\top
P^\varepsilon(t)\gamma_{x}^{\bm{p}}(y^\varepsilon(t),u^\varepsilon(t))
+\sum_{d=1}^{\infty}\sum\limits_{\bm{p}\in\mathbb{N}_d^n}\gamma_{x}^{\bm{p}}(y^\varepsilon(t),u^\varepsilon(t))^\top R^{\bm{p},\varepsilon(t)}\nonumber\\
&&\left.+\sum_{d=1}^{\infty}\sum\limits_{\bm{p}\in\mathbb{N}_d^n}R^{\bm{p},\varepsilon}(t)\gamma_{x}^{\bm{p}}(y^\varepsilon(t),u^\varepsilon(t))
+H_{xx}(y^\varepsilon(t),u^\varepsilon(t),\lambda^\varepsilon,p^\varepsilon(t),J^\varepsilon(t))\right]dt\nonumber\\
&&-\sum_{d=1}^{\infty}\sum\limits_{\bm{p}\in\mathbb{N}_d^n}R^{\bm{p},\varepsilon}(t)dH^{\bm{p}}(t)\nonumber\\
P^\varepsilon(T)&=&\lambda^\varepsilon h_{xx}(y_0^\varepsilon,y^\varepsilon(T))
 +\sum_{j=1}^k\mu^{\varepsilon j}G_{xx}^j(y_0^\varepsilon,y^\varepsilon(T)).
\end{eqnarray}
have unique solutions
$(p^\varepsilon(\cdot),\{J^{\bm{p},\varepsilon}(\cdot)\}_{\bm{p}\in\mathbb{N}^n})$
and $(P^\varepsilon(\cdot),\{R^{\bm{p},\varepsilon}(\cdot)\}_{\bm{p}\in\mathbb{N}^n})$
respectively, with $p^\varepsilon(\cdot)$ and $P^\varepsilon(\cdot)$
being cadlag processes.

Using It$\hat{o}$'s formula, we have from (4.22), (4.28) and (4.29), that
\begin{eqnarray}
\begin{array}{cl}
&E<\lambda^\varepsilon
h_x(y_0^\varepsilon,y^\varepsilon(T))+\sum\limits_{j=1}^k\mu^{\varepsilon
j}G_x^j(y_0^\varepsilon,y^\varepsilon(T))+\int_0^T\lambda^\varepsilon\ell_x(y^\varepsilon(s),u^\varepsilon(s)),
y_1^\varepsilon(T)+y_2^\varepsilon(T)>\\
=&E<p^\varepsilon(T),y_1^\varepsilon(T)+y_2^\varepsilon(T)>\\
=&<p^\varepsilon(0),\eta>|I_\rho|+E\int_0^T(p(s),\triangle g^\varepsilon(s,u^{\varepsilon\rho}(s)))ds\\
&+\sum\limits_{d=1}^\infty\sum\limits_{\bm{p}\in\mathbb{N}_d^n}E\int_0^T(J^{\bm{p}}(s),
\triangle\gamma^{\varepsilon,\bm{p}}(s,u^{\varepsilon\rho}(s)))ds\\
&+\frac{1}{2}E\int_0^T(p(s),g_{xx}(y^\varepsilon(s),u^\varepsilon(s))y_1^\varepsilon(s)y_1^\varepsilon(s))ds\\
&+\frac{1}{2}\sum\limits_{d=1}^\infty\sum\limits_{\bm{p}\in\mathbb{N}_d^n}E\int_0^T(J^{\bm{p}}(s),
\gamma^{\bm{p}}(y^\varepsilon(s),u^\varepsilon(s))y_1^\varepsilon(s)y_1^\varepsilon(s))ds\\
&+\sum\limits_{d=1}^\infty\sum\limits_{\bm{p}\in\mathbb{N}_d^n}E\int_0^T(J^{\bm{p}}(s),
\triangle\gamma^{\varepsilon,\bm{p}}_x(s,u^{\varepsilon\rho}(s))y_1^\varepsilon(s))ds
\end{array}
\end{eqnarray}

Applying Ito's formula to the matrix-valued processes
\begin{eqnarray}
Y(s)=y_1(s)y_1^\top(s)=\left(
                      \begin{array}{ccc}
                        y_1^1y_1^1 & \ldots & y_1^1y_1^m \\
                        \vdots & \vdots& \vdots \\
                        y_1^1y_1^m & \ldots & y_1^my_1^m \\
                      \end{array}
                    \right)\nonumber
\end{eqnarray}
we have
\begin{eqnarray}
dY(t)&=&\left[Y(t)g_x^\top(t)+g_x(t)Y(t)
+\sum_{d=1}^{\infty}\sum\limits_{\bm{p}\in\mathbb{N}_d^n}\gamma_{x}^{\bm{p}}(t)Y(t)\gamma_{x}^{\bm{p}}(t)^\top+\Phi^\varepsilon(t)\right]dt\nonumber\\
&&+\sum_{d=1}^{\infty}\sum\limits_{\bm{p}\in\mathbb{N}_d^n}\left[Y(t)\gamma_x^{\bm{p}}(t)^\top+\gamma_x^{\bm{p}}(t)Y(t)
+\Omega^{\bm{p},\varepsilon}(t)\right] dH^{\bm{p}}(t)
\end{eqnarray}
where
\begin{eqnarray}
\Phi^\varepsilon(t)&=&\sum_{d=1}^{\infty}\sum\limits_{\bm{p}\in\mathbb{N}_d^n}\gamma_x^{\bm{p}}(t)y_1(t)\triangle\gamma^{\bm{p}}(t,u^\varepsilon(t))
^\top+\sum_{d=1}^{\infty}\sum\limits_{\bm{p}\in\mathbb{N}_d^n}\triangle\gamma^{\bm{p}}(t,u^\varepsilon(t))y_1(t)^T
\gamma_x^{\bm{p}}(t)^\top\nonumber\\
&&+\sum_{d=1}^{\infty}\sum\limits_{\bm{p}\in\mathbb{N}_d^n}\triangle\gamma^{\bm{p}}(t,u^\varepsilon(t))
\sum_{d=1}^{\infty}\sum\limits_{\bm{p}\in\mathbb{N}_d^n}\triangle\gamma^{\bm{p}}(t,u^\varepsilon(t))^\top\nonumber\\
\Omega^{\bm{p},\varepsilon}(t)&=&y_1(t)\triangle\gamma^{\bm{p}}(t,u^\varepsilon(t))^\top
+\triangle\gamma^{\bm{p}}(t,u^\varepsilon(t))y_1(t)^\top\nonumber\\
&&+\sum_{d=1}^{\infty}\sum\limits_{\bm{p}\in\mathbb{N}_d^n}\triangle\gamma^{\bm{p}}(t,u^\varepsilon(t))
\sum_{d=1}^{\infty}\sum\limits_{\bm{p}\in\mathbb{N}_d^n}\triangle\gamma^{\bm{p}}(t,u^\varepsilon(t))
^\top\nonumber
\end{eqnarray}
and
\begin{eqnarray}
\begin{array}{cl}
&\lambda^\varepsilon
\mathbb{E}y_1^{\varepsilon,\top}(T)h_{xx}(y_0^\varepsilon,y^\varepsilon(T))y_1^\varepsilon(T)
 +\sum\limits_{j=1}^k\mu^{\varepsilon
 j}\mathbb{E}y_1^{\varepsilon,\top}(T)G_{xx}^j(y_0^\varepsilon,y^\varepsilon(T))y_1^\varepsilon(T)\\
=&trE[P^\varepsilon(T)y_1^\varepsilon(T)y_1^{\varepsilon,\top}(T)]\\
=&-E\int_0^Ty_1^{\varepsilon,\top}(s)H_{xx}(y^\varepsilon(s),u^\varepsilon(s),\lambda^\varepsilon,
p^\varepsilon(s),J^\varepsilon(s))y_1^\varepsilon(s)ds\\
&+E\int_0^Ttr
P^\varepsilon(s)\left[\sum\limits_{d=1}^\infty\sum\limits_{\bm{p}\in\mathbb{N}_d^n}
\triangle\gamma^{\varepsilon,\bm{p}}(s;u^{\varepsilon\rho}(s))\sum\limits_{d=1}^\infty\sum\limits_{\bm{p}\in\mathbb{N}_d^n}
\triangle\gamma^{\varepsilon,\bm{p},\top}(s;u^{\varepsilon\rho}(s))\right]ds\\
&+E\int_0^Ttr\left[ \sum\limits_{d=1}^\infty\sum\limits_{\bm{p}\in\mathbb{N}_d^n}
R^{\bm{p},\varepsilon}(s)\right]^\top\left[\sum\limits_{d=1}^\infty\sum\limits_{\bm{p}\in\mathbb{N}_d^n}
\triangle\gamma^{\varepsilon,\bm{p}}(s;u^{\varepsilon\rho}(s))\sum\limits_{d=1}^\infty\sum\limits_{\bm{p}\in\mathbb{N}_d^n}
\triangle\gamma^{\varepsilon,\bm{p},\top}(s;u^{\varepsilon\rho}(s))\right]ds\\
&+2E\int_0^Ttr
P^\varepsilon(s)\left[\sum\limits_{d=1}^\infty\sum\limits_{\bm{p}\in\mathbb{N}_d^n}\gamma_x^{\varepsilon,\bm{p}}(s)y_1^\varepsilon(s)
\triangle\gamma^{\varepsilon,\bm{p},\top}(s;u^{\varepsilon\rho}(s))\right]ds
\end{array}
\end{eqnarray}

Noting the estimates (4.24), we conclude from (4.25)-(4.27) and (4.30)-(4.32) that
\begin{eqnarray}
\begin{array}{cl}
&E<\lambda^\varepsilon
h_y(y_0^\varepsilon,y^\varepsilon(T))+\sum\limits_{j=1}^k\mu^{\varepsilon
j}G_y^j(y_0^\varepsilon,y^\varepsilon(T)+p^\varepsilon(0),\eta>|I_\rho|\\
&+\int_0^Tl^\varepsilon(s,u^{\varepsilon\rho})ds+o(|I_\rho|)\geq
-|I_\rho|\sqrt{\varepsilon+2\delta(\varepsilon)}\sqrt{1+|\eta|^2}
\end{array}
\end{eqnarray}
where $l^\varepsilon(\cdot;v)$ is defined by
\begin{eqnarray}
\begin{array}{ccl}
l^\varepsilon(s;v)&=:&E\left(H(y(s),u^\varepsilon(s),\lambda,p(s),J(s))-H(y(s),u(s),\lambda,p(s),J(s))\right)\\
&&+\frac{1}{2}Etr
P^\varepsilon(s)\left[\sum\limits_{d=1}^\infty\sum\limits_{\bm{p}\in\mathbb{N}_d^n}
\triangle\gamma^{\varepsilon,\bm{p}}(s;u^{\varepsilon\rho}(s))\sum\limits_{d=1}^\infty\sum\limits_{\bm{p}\in\mathbb{N}_d^n}
\triangle\gamma^{\varepsilon,\bm{p},\top}(s;u^{\varepsilon\rho}(s))\right]\\
&&+\frac{1}{2}Etr\left[ \sum\limits_{d=1}^\infty\sum\limits_{\bm{p}\in\mathbb{N}_d^n}
R^{\varepsilon,\bm{p}}(s)\right]^\top\left[\sum\limits_{d=1}^\infty\sum\limits_{\bm{p}\in\mathbb{N}_d^n}
\triangle\gamma^{\varepsilon,\bm{p}}(s;u^{\varepsilon\rho}(s))\sum\limits_{d=1}^\infty\sum\limits_{\bm{p}\in\mathbb{N}_d^n}
\triangle\gamma^{\varepsilon,\bm{p},\top}(s;u^{\varepsilon\rho}(s))\right]
\end{array}
\end{eqnarray}

\textsl{Step 3. Differentiability.} For given $v(\cdot)\in \mathcal{U}_{ad}$, applying Lemma 3.2 to the real valued Lebesgue integrable function, we know that there exists $I_\rho\subset[0,T]$ such that
\begin{eqnarray}
\begin{array}{rcl}
  |I_\rho| &=&\rho , \\
  \int_{I_\rho}\ell^\varepsilon(s,v(s))ds &=&\rho\int_0^T\ell^\varepsilon(s;v(s))ds+o(\rho),\quad as \quad \rho\rightarrow 0
\end{array}
\end{eqnarray}
Next choose the above $I_\rho$ in (4.20), and we have
\begin{eqnarray}
 \int_{I_\rho}\ell^\varepsilon(s,v(s))ds &=&\rho\int_0^T\ell^\varepsilon(s;u^{\varepsilon\rho}(s))ds
\end{eqnarray}
From (4.33)-(4.36), we conclude for given $v(\cdot)\in\mathcal{U}_{ad}$ that
\begin{eqnarray}
\begin{array}{l}
E<\lambda^\varepsilon
h_y(y_0^\varepsilon,y^\varepsilon(T))+\sum\limits_{j=1}^k\mu^{\varepsilon,
j}G_y^j(y_0^\varepsilon,y^\varepsilon(T)+p^\varepsilon(0),\eta>\rho+\rho\int_0^T\ell^\varepsilon(s,v(s))ds\\
\geq-\rho\sqrt{\varepsilon+2\delta(\varepsilon)}\sqrt{1+|\eta|^2}+o(\rho),\quad as \quad \rho\rightarrow 0.
\end{array}
\end{eqnarray}
Hence
\begin{eqnarray}
\begin{array}{l}
E<\lambda^\varepsilon
h_y(y_0^\varepsilon,y^\varepsilon(T))+\sum\limits_{j=1}^k\mu^{\varepsilon,
j}G_y^j(y_0^\varepsilon,y^\varepsilon(T)+p^\varepsilon(0),\eta>+\int_0^T\ell^\varepsilon(s,v(s))ds\\
\geq-\sqrt{\varepsilon+2\delta(\varepsilon)}\sqrt{1+|\eta|^2},\quad \forall \eta\in\mathbb{R}^m \quad \forall v(\cdot)\in\mathcal{U}_{ad}.
\end{array}
\end{eqnarray}
This implies that
\begin{eqnarray}
\begin{array}{l}
\lambda^\varepsilon
Eh_y(y_0^\varepsilon,y^\varepsilon(T))+\sum\limits_{j=1}^k\mu^{\varepsilon,
j}EG_y^j(y_0^\varepsilon,y^\varepsilon(T)+p^\varepsilon(0)\leq C\sqrt{3\varepsilon},\\
\int_0^T\ell^\varepsilon(s,v(s))ds\geq-\sqrt{\varepsilon+2\delta(\varepsilon)},\quad \forall v(\cdot)\in\mathcal{U}_{ad}.
\end{array}
\end{eqnarray}

\textsl{Step 4. Passing to the limit.} Without loss of generality,
we assume that $\lambda^\varepsilon\rightarrow\lambda$,
$\mu^\varepsilon\rightarrow\mu$, as $\varepsilon\rightarrow 0+$.

Let $\varepsilon\rightarrow 0+$. Equation (4.19)$_2$ gives the following:
\begin{eqnarray}
\begin{array}{l}
E\int_0^T\left(H(y(s),u^\varepsilon(s),\lambda,p(s),J(s))-H(y(s),u(s),\lambda,p(s),J(s))\right)ds\\
+\frac{1}{2}E\int_0^Ttr
P^\varepsilon(s)\left[\sum\limits_{d=1}^\infty\sum\limits_{\bm{p}\in\mathbb{N}_d^n}
\triangle\gamma^{\varepsilon,\bm{p}}(s;u^{\varepsilon\rho}(s))\sum\limits_{d=1}^\infty\sum\limits_{\bm{p}\in\mathbb{N}_d^n}
\triangle\gamma^{\varepsilon,\bm{p},\top}(s;u^{\varepsilon\rho}(s))\right]ds\\
+\frac{1}{2}E\int_0^Ttr\left[ \sum\limits_{d=1}^\infty\sum\limits_{\bm{p}\in\mathbb{N}_d^n}
R^{\varepsilon,\bm{p}}(s)\right]^\top\left[\sum\limits_{d=1}^\infty\sum\limits_{\bm{p}\in\mathbb{N}_d^n}
\triangle\gamma^{\varepsilon,\bm{p}}(s;u^{\varepsilon\rho}(s))\sum\limits_{d=1}^\infty\sum\limits_{\bm{p}\in\mathbb{N}_d^n}
\triangle\gamma^{\varepsilon,\bm{p},\top}(s;u^{\varepsilon\rho}(s))\right]ds\\
\geq 0,\quad\forall v(\cdot)\in U_{ad};
\end{array}
\end{eqnarray}
this implies (4.11). Furthermore, (4.11) is obtained from (4.19)$_1$,(4.9)$_2$ is obtained from (4.39)$_1$, and the rest of Theorem 4.1 is checked from (4.28) and (4.29).

\textsl{Step 5. The unbounded case of $\mathcal{U}_{ad}$ in $L_{\mathscr{F},p}^{\infty,8}[[0,T];
\mathbb{R}^m]$.}The proof procedure is the same as the step 5 in Tang and Li [34].

The proof of Theorem 4.1 is complete.$\Box$

\section{Conclusions}

In this paper, necessary maximum principle for optimal control of stochastic system driven by multidimensional Teugel's martingales is proved, where the multidimensional Teugel's martingales are constructed by orthogonalizing the multidimensional L\'{e}vy processes. The control variable is allowed to enter the coefficients of the Teugel's martingales, and the control domain is nonconcave. The technique for proving the maximum principle and the obtained result are almost similar to Peng [26] and Tang and Li [34].

\

\end{document}